\documentclass[11pt]{amsart}
\usepackage{a4wide,amsmath,amssymb,amscd,verbatim,
psfrag,subfigure,graphicx}

\newcommand{\lebn}

\newtheorem{theorem}[equation]{Theorem}
\newtheorem{corollary}[equation]{Corollary}

\theoremstyle{definition}

\newtheorem{definition}[equation]{Definition}

\newtheorem{example}[equation]{Example}

\numberwithin{equation}{section}

\newcommand{\C}{\mathbb{C}}

\newcommand{\R}{\mathbb{R}}

\newcommand{\Z}{\mathbb{Z}}

\begin{document}
\bibliographystyle{plain}
\title[Some examples in toric geometry \lebn]
{Some examples in toric geometry}
\author{Yusuf Civan}
\address{Department of Mathematics, Suleyman Demirel University,
Isparta, 32260, Turkey}   
\email{ycivan@fef.sdu.edu.tr}



\begin{abstract}
We present various examples in toric geometry concerning the relationship between 
smooth toric varieties and quasitoric manifolds (or more generally unitary 
torus manifolds), and extend the results of 
\cite{GKS} to prove the non-existence of almost complex quasitorics over the 
duals of some certain cyclic $4$-polytopes. We also provide
the sufficient conditions on the base polytope and the characteristic 
map so that the resulting quasitoric manifold is almost complex, answering the 
question proposed by Davis $\&$ Januszkiewicz \cite{DJ} .   
\end{abstract}

\maketitle
\section{Introduction}\label{intro}
There is an unfortunate clash on the notion of a \emph{toric manifold}; for instance, a 
toric manifold for geometers is a smooth toric variety, however, for an algebraic 
topologist, it is a smooth, even dimensional real manifold acted upon by a compact
torus such that the quotient is homeomorphic to a simple convex polytope as a manifold
with corners. At first, the class of the latter objects seems to contain the formers
(indeed,  Davis $\&$ Januszkiewicz \cite{DJ} claimed to do so), we provide an example
showing that (see Example \ref{exa2}) it is not the case. 
On the other hand, it is known that any smooth toric 
variety associated with a fan arising from a polytopal simplicial complex is a toric 
manifold in the latter sense \cite{YC}. Conversely, our Example \ref{exa1} exhibits
an almost complex quasitoric manifold that does not arise from a smooth toric variety. 
To distinguish these two classes, Buchstaber $\&$ Panov \cite{BT} prefer to call a 
toric manifold of Davis $\&$ Januszkiewicz a \emph{quasitoric manifold}, so do we, 
and preserve the term \emph{toric manifold} for smooth toric varieties.  

The geometric and computational flavor of the theory enables us to translate 
topological problems into combinatorial ones and vice verse. In this guise, the 
existence of an almost complex structure compatible with the action on a given 
quasitoric manifold may be formulated as a combinatorial property carried by the 
associated simple polytope and the characteristic map, which we characterize 
in Theorem \ref{T1} .   
 
In recent years, Masuda's work on unitary torus manifolds \cite{MM} contributes to the
theory from an unfashionable manner. Instead of starting with some combinatorial 
ingredients, he begins with a closed, connected, stably complex manifold $M^{2n}$
equipped with a $T^n$-action such that the $T^n$-fixed point set is nonempty and
isolated, and then he recovers a combinatorial object called \emph{multi-fan}
associated to $M^{2n}$. It can be verified that the class of unitary torus 
manifolds contains all quasitoric and toric manifolds. However, this containment
is strict by Example \ref{exa4}. Moreover, under certain restrictions, 
we can easily adopt our programme in order to characterize the existence of 
an almost complex structure on unitary torus manifolds (see Corollary \ref{C1}). 

\section{Almost complex quasitorics}\label{Acq}

Our first purpose here is to characterize the existence of an almost complex structure
on quasitoric manifolds in terms of the related simple polytopes and dicharacteristic maps.
So we begin with introducing some notations, and we refer readers to \cite{BT} for a more 
detailed expositions in the theory.  

Let $M^{2n}$ be a quasitoric manifold over $P^n$ and let
$\mathcal{F}=\{F_1,\ldots, F_m\}$ be the set of facets of $P^n$.
Then for each $F_i$, the pre-image $\pi^{-1}(F_i)$ is a submanifold
$M_i^{2(n-1)}\subset M^{2n}$ with isotropy group a circle $T(F_i)$ in $T^n$.
Since there is a one-to-one correspondence (up to a sign) between the set of
primitive vectors in $\Z^n$ and the subcircles in $T^n$, we obtain the 
\emph{characteristic map} of  $M^{2n}$ given by 
\begin{align*}
\mathbf{\lambda}\colon& \mathcal{F}\rightarrow \Z^n\\
&F_i\mapsto \mathbf{\lambda}(F_i):=\mathbf{\lambda}_i,
\end{align*}
where $\lambda_i$ generates the circle $T(F_i)$ in $T^n$.
We note that the map $\mathbf{\lambda}$ is well defined only up to a sign, and 
if the sign of each $\mathbf{\lambda}_i$ is chosen, we then call
$\mathbf{\lambda}$ a \emph{dicharacteristic map} of $M^{2n}$. Therefore, there are 
$2^m$ dicharacteristic maps in total attached to $M^{2n}$. On the other hand,
each such choice for $\mathbf{\lambda}_i$ determines an orientation of the normal
bundle $\nu_i$ of $M_i^{2(n-1)}$, so an orientation for $M_i^{2(n-1)}$. Conversely,
an \emph{omniorientation} of $M^{2n}$ consists of a choice of an orientation 
for every submanifold $M_i^{2(n-1)}$, which in turn settles a sign for each vector
$\mathbf{\lambda}_i$. Thus, every omniorientation is equipped with a unique 
dicharacteristic map and vice versa. Buchstaber $\&$ Ray \cite{BR} were able to 
show that any omniorientation of $M^{2n}$ induces a stably complex structure on it
by means of the following isomorphism:

\begin{equation}
\tau(M^{2n})\oplus \R^{2(m-n)}\cong \rho_1\oplus \ldots \oplus \rho_m,
\end{equation}   
where $\rho_i$ is the pull back of the line bundle corresponding to the Thom class 
defined by $\nu_i$ along the Pontryagin-Thom collapse.

Since $P^n$ is simple, each vertex $v$ of $P^n$ can be written
as an intersection of $n$ facets:
\begin{equation}\label{eq1}
v=F_{i_1}\cap \ldots \cap F_{i_n}.
\end{equation}

Assign to each facet $F_{i_k}$ the edge $E_k:\bigcap_{j\neq k} F_{i_j}$, and let
$\mathbf{e}_k$ be the vector along $E_k$ beginning at $v$. Then, depending on the
ordering of the facets \eqref{eq1}, the vectors $\mathbf{e}_1, \ldots, \mathbf{e}_n$
form either positively or negatively oriented basis of $\R^n$. Throughout 
this ordering is assumed to be so that $\mathbf{e}_1, \ldots, \mathbf{e}_n$ is a
positively oriented basis.

For a given dicharacteristic map $\mathbf{\lambda}$, we define $\Lambda$ to be the 
$(n\times m)$-matrix whose $i$-th column is formed by the vector $\mathbf{\lambda}_i^t$
for any $1\leq i\leq m$. We let $\Lambda_v:=\Lambda_{i_1,\ldots, i_n}$ denote 
the maximal minor of $\Lambda$ formed by the columns $i_1,\ldots, i_n$, where 
$v=F_{i_1}\cap \ldots \cap F_{i_n}$. From the definition of a 
characteristic map, we have that 
\begin{equation*}
\mathrm{det}\;\Lambda_v=\mp 1
\end{equation*}
for any vertex $v\in P^n$.

\begin{definition}
The \emph{sign} of a vertex $v\in P^n$ is defined to be
\begin{equation*}
\sigma(v):=\mathrm{det}\;\Lambda_v
\end{equation*}
\end{definition}

\begin{theorem}\label{T1} 
An omniorientation of a quasitoric manifold $M^{2n}$ over $P^n$  
arises from a $T^n$-invariant almost complex structure on $M^{2n}$ if and only if 
$\sigma(v)=1$ for each vertex $v\in P$.
\end{theorem}
\begin{proof}
The necessity part of the claim has already appeared in \cite{BT}, so for the sufficiency,
assume that $\sigma(v)=1$ for each vertex $v\in P$. However, this guaranties that 
the Euler number of the resulting stably complex structure equals to 
the Euler number of $M^{2n}$; hence, it arises from an almost complex structure 
by the Proposition $4.1$ of \cite{AT}. It means that the complex structure $J$ on
$\tau(M^{2n})\oplus \R^{2(m-n)}$ splits as $J=(J_1, J_2)$, where $J_1$ and  $J_2$
are complex structures on $\tau(M^{2n})$ and the trivial portion respectively.
However, since $J$ is $T^n$-invariant, so is $J_1$. 

On the other hand, the assumption
$\sigma(v)=1$ for each vertex $v\in P$ prevails the fact that each $T(F_i)$-fixed 
submanifolds $M_i^{2(n-1)}$ has an almost complex structure induced by that of $M^{2n}$. 
This may be achieved systematically  by obtaining the dicharacteristic map of 
$M_i^{2(n-1)}$ from that of $M^{2n}$ such a way that for each fixed point of $M_i^{2(n-1)}$, 
the sign of the corresponding vertex is equal to $1$.
\end{proof}
\section{Examples}\label{Exa}

Once we described almost complex structures on quasitorics, it would be of interest
to find an example of an almost complex quasitoric that is not a toric manifold, 
since all known almost complex quasitoric manifolds at the moment are also toric 
manifolds (see Problem $2.2.11$ of \cite{BT}).  

\begin{figure}[ht]
\begin{center}
\psfrag{a}{$v_1$}
\psfrag{b}{$v_2$}
\psfrag{c}{$v_3$}
\psfrag{d}{$v_4$}
\psfrag{e}{$v_5$}
\psfrag{f}{$(1,1)$}
\psfrag{A}{$\lambda_1=(0,-1)$}
\psfrag{B}{$\lambda_2=(1,1)$}
\psfrag{C}{$\lambda_3=(1,2)$}
\psfrag{D}{$\lambda_4=(-2,-3)$}
\psfrag{E}{$\lambda_5=(-1,-2)$}
\includegraphics[width=3.4in,height=2.7in]{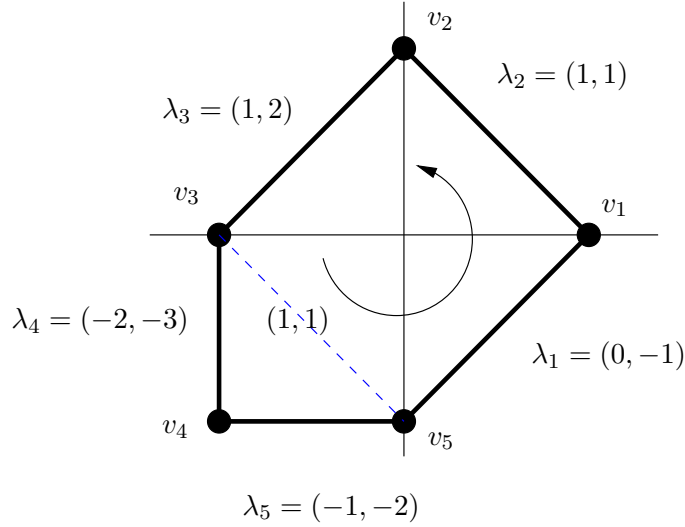}
\end{center}
\caption{The base polytope and a dicharacteristic map for 
$\C P^2 \# \C P^2 \# \overline{\C P}^2$}\label{AQ}
\end{figure}

\begin{example}\label{exa1}
Let $M^4=\C P^2 \# \C P^2 \# \overline{\C P}^2$ be the quasitoric $4$-manifold over the
pentagon $P$ oriented counterclockwise shown by the Figure \ref{AQ} with the given 
dicharacteristic map (compare to the list given in [\cite{OR}, p.$552$], where $\#$ denotes
the connected sum. The triangle part of $P$ corresponds to the base polytope of $\C P^2$, 
and the quadruple portion for $\C P^2 \# \overline{\C P}^2$, where $\overline{\C P}^2$ 
denotes the complex projective plane with the reversed orientation. 

It is easy to check that $\sigma(v_i)=1$ for any $1\leq i \leq 5$. Hence, by the Theorem
\ref{T1}, the resulting stably complex structure on $M^4$ is induced by a $T^2$-invariant 
almost complex structure. However, $M^4$ can not have the diffeomorphic type (or even 
homeomorphic type) of a toric manifold by the classification theorem of
Fischli $\&$ Yavin \cite{FY}. Therefore,  $M^4$ is an almost complex quasitoric manifold
which is not a toric manifold.    
\end{example}

We next present a toric manifold that is not a quasitoric manifold.

\begin{example}\label{exa2}
Let $\mathcal{B}^3$ denote the Barnette sphere, which is a star-shaped, non-polytopal 
simplicial $3$-sphere (see \cite{DB} or \cite{GE}). Since it is star-shaped, it spans a
complete fan $\Sigma(\mathcal{B})$ in $\R^4$ that may not be necessarily smooth. However,
applying a finite number of stellar subdivision, we can turn $\Sigma(\mathcal{B})$ into a 
smooth fan, which we still continue to denote by $\Sigma(\mathcal{B})$. Since the 
non-polytopality of $\mathcal{B}^3$ comes from the fact that it has a double edge,
the stellar subdivision on $\mathcal{B}^3$ will not remove such an edge so that the 
underlying simplicial sphere of $\Sigma(\mathcal{B})$ is still non-polytopal. If we denote
the associated toric manifold by $X(\mathcal{B})$, it contains an algebraic torus as a
dense subset that acts on $X(\mathcal{B})$ smoothly. Furthermore, the quotient of 
$X(\mathcal{B})$ by the action of the compact $4$-torus included in the algebraic torus 
is a $4$-ball whose facial structure is not isomorphic to a simple convex polytope. Thus, 
the toric manifold $X(\mathcal{B})$ can not be a quasitoric manifold.     
\end{example}

\section{Quasitorics and cyclic polytopes}\label{Qcp}
During the $1990$s much of the effort in toric geometry was spent on the classification
and the existence of toric manifolds in a purely combinatorial sense. From the 
geometric point of view, for a given simplicial complex $\Delta$, the existence of a 
fan $\Sigma(\Delta)$ associated with it depends on the star-shapness of $\Delta$. However,
the smoothness of such a fan can be expressed as a realizability problem in synthetic 
geometry. In other words, we have to have a realization of $\Sigma(\Delta)$ in such a way 
that some determinant conditions are satisfied. In the projective case, where the complex
$\Delta$ is the boundary complex of a simplicial convex polytope, the neighborliness of the
polytope puts more restriction on the existence of such realization. In fact, this is
the key point on which Gretenkort, et al \cite{GKS} was able to prove that there does
not exits a smooth fan whose spanning simplicial complex is the boundary complex of a cyclic
polytope with $n\geq 7$ vertices. Even though, the basic combinatorial ingredients for 
constructing quasitorics seem to be more flexible, when we require almost complex 
structures on them, similar synthetic geometry problems will appear. 

Instead of using the moment curve to construct a cyclic $4$-polytope $C_4(n)$ with 
$n$ vertices, we may alternatively use the \emph{Carath$\acute{e}$odory curve} 
(see \cite{GZ});
\begin{align*}
\mathbf{p}\colon & \R \rightarrow \R^4 \\
&u\mapsto \mathbf{p}(u):=(cos(u), sin(u), cos(2u), sin(2u))
\end{align*}
so that $C_4(n)=\mathrm{conv}\{\mathbf{p}(t_1),\ldots, \mathbf{p}(t_n)\}$
for $0\leq t_1< t_2<\ldots < t_n< 2\pi$. We then denote by $D_4(n)$
the \emph{dual} (or the \emph{polar} when $0\in \mathrm{relint}(C_4(n))$) of $C_4(n)$.
It follows that $D_4(n)$ is a simple convex $4$-polytope.

\begin{example}\label{exa3}
Consider 
$0< \frac{\pi}{4}<\frac{\pi}{2}<\frac{3\pi}{4}<\pi<
\frac{5\pi}{4}<\frac{3\pi}{2}< 2\pi$, and let $C_4(7)$ be the cyclic
polytope with vertices
\begin{align*}
&v_1=\mathbf{p}(0),\; v_2=\mathbf{p}(\frac{\pi}{4}),\; v_3=\mathbf{p}(\frac{\pi}{2}),\;
v_4=\mathbf{p}(\frac{3\pi}{4}),\\
&v_5=\mathbf{p}(\pi),\;v_6=\mathbf{p}(\frac{5\pi}{4}),\;v_7=\mathbf{p}(\frac{3\pi}{2}).
\end{align*}
It is easy to verify that $0\in \mathrm{relint}(C_4(7))$ so that its polar $D_4(7)$ exists.
By the anti-isomorphism between the face lattices of $C_4(7)$ and $D_4(7)$, we denote the 
facets of $D_4(7)$ by $F_1, F_2, F_3, F_4, F_5, F_6, F_7$ corresponding to the vertices
of $C_4(7)$ with the same indexes. By the Gale's evenness condition and the choice of our 
realization, the list of positively ordered facets meeting at some vertex of $D_4(7)$ 
can be given as follows:
\begin{align*}
&1234 &1267  &&2345  &&&3467 \\
&2137 &3147  &&2356  &&&4567 \\
&2145 &4157  &&2367       \\
&1256 &1567  &&3456       \\     
\end{align*}

By analogy with the proof given in [\cite{GKS}, p.257], it can be shown that for any 
dicharacteristic map on $D_4(7)$, there is at least one vertex $v$ in the above list 
such that $\sigma(v)=-1$. Therefore, there does not exist a dicharacteristic
map on $D_4(7)$ satisfying $\sigma(v)=1$ for each vertex $v\in D_4(7)$; hence,
there is no almost complex quasitoric manifold with the base polytope $D_4(7)$
realized as above. However, we may construct a stably complex quasitoric over $D_4(7)$
with the dicharacteristic map given, for example by
\begin{align*}
&\lambda(F_1)=(0,1,0,0),\; \lambda(F_2)=(1,0,0,0),\; \lambda(F_3)=(0,0,1,0),\; 
\lambda(F_4)=(-1, 0, -1,-1),\\
&\lambda(F_5)=(1,-1,0,-1),\; \lambda(F_6)=(1,-1,-1,0),\; \lambda(F_7)=(0,0,0,1).
\end{align*}
\end{example}

\begin{theorem}\label{T2}
There does not exist an almost complex quasitoric manifold over the polytope
$D_4(n)$ with $n\geq 7$.  
\end{theorem}
\begin{proof}
We first note that the existence of an almost complex structure is independent of
any specific geometric realization of the polytope. Therefore, when $n=7$, 
the claim follows from the Example \ref{exa3}. A similar reason applies 
to the case $n\geq 8$.  
\end{proof} 

\section{Unitary torus manifolds}\label{Utm}
Since the class of multi-fans contains all convex polytopes as well as fans,
the unitary torus manifolds may be thought of a generalization of quasitoric and toric
manifolds. In this circumstance, it would be interesting to clarify this containment.
In other words, it is not clear whether there exists a unitary torus manifold
that is neither a quasitoric nor a toric manifold. We provide such an example in
Example \ref{exa4}. We may also extend our discussion in Section \ref{Acq} to 
deal with the existence of an almost complex structure on unitary torus manifolds 
under some restrictions. We refer readers to \cite{HM} for a more detailed exposition 
and notation on unitary torus manifolds.     

\begin{example}\label{exa4}
Let  $\mathcal{B}^3$ denote the Barnette sphere with the following $3$-simplicies (see \cite{GE}):
\begin{equation*}
\begin{split}
&[x_1, x_2, x_3, x_4]\quad [x_3, x_4, x_5, x_6]\quad [x_1, x_2, x_5, x_6]\quad  [x_1, x_2, x_4, x_7]
\quad [x_1, x_3, x_4, x_7] \quad [x_3, x_4, x_6, x_7]\\ 
&[x_3, x_5, x_6, x_7]\quad  [x_1, x_2, x_5, x_7]\quad [x_2, x_5, x_6, x_7]\quad [x_2, x_4, x_6, x_7]
\quad [x_1, x_2, x_3, x_8]\quad [x_2, x_3, x_4, x_8]\\
&[x_3, x_4, x_5, x_8]\quad [x_4, x_5, x_6, x_8]\quad [x_1, x_2, x_6, x_8]\quad [x_1, x_5, x_6, x_8]
\quad [x_1, x_3, x_5, x_8]\quad [x_2, x_4, x_6, x_8]\\
&\mathrm{and}\quad [x_1, x_3, x_5, x_7]\quad \mathrm{as\; the \;base}.
\end{split} 
\end{equation*}
The $f$ and $h$-vectors of $\mathcal{B}$ are given by $f(\mathcal{B})=(8, 27, 38, 19)$ and
$h(\mathcal{B})=(1, 4, 9, 4, 1)$. We define 
$\lambda \colon V(\mathcal{B})\rightarrow \Z^4$ by
\begin{equation*}
\begin{split}
&\lambda(x_1)=(1, 0, 0, 0),\quad \lambda(x_2)=(0, 1, -1, 2),\quad \lambda(x_3)=(0, 1, 0, 0),\\
&\lambda(x_4)=(0, 0, 1, -1),\quad \lambda(x_5)=(0, 0, 1, 0),\quad  \lambda(x_6)=(1, -1, 0, -1),\\
&\lambda(x_7)=(0, 0, 0, 1),\quad \mathrm{and}\quad \lambda(x_8)=(1, 0, 0, -1),
\end{split}
\end{equation*}
and denote by $\Lambda$, the matrix whose columns consist of vectors $\lambda(x_i)^t$ for 
$1\leq i \leq 8$. It then follows that 
\begin{equation*}
\mathrm{det}\;\Lambda_{\sigma}=\mp 1,
\end{equation*}
for each $3$-simplex $\sigma$ of $\mathcal{B}$, where $\Lambda_{\sigma}$ is the maximal minor of
$\Lambda$ corresponding to the simplex $\sigma$. Since $\mathcal{B}$ is a simplicial sphere, 
the pair $(\mathcal{B}, \lambda)$ defines a unitary torus manifold $M_{\lambda}(\mathcal{B})$.
It is obvious that $\mathcal{B}$ can not span a smooth fan in $\R^4$ with the generating set 
$\{\lambda(x_i)\colon  1\leq i \leq 8\}$ so that $M_{\lambda}(\mathcal{B})$ can not be a toric manifold. 
Similarly, since $\mathcal{B}$ is a non-polytopal sphere; the manifold $M_{\lambda}(\mathcal{B})$ 
is not quasitoric.
\end{example}
We note that any unitary torus $4$-manifold is quasitoric by the Steinitz Theorem which
asserts that any simplicial $2$-sphere is polytopal.

We next determine the condition for a unitary torus manifold in order to carry
an almost complex structure under some restrictions. Our limitation here only requires 
that the Euler characteristic of the manifold equals to the sum of the $h$-vector of the 
underlying multi-fan. 

Let $M^{2n}$ be a unitary torus manifold such that $M_{I(p)}=p$ for any $p\in M^T$, and 
assume that if $p=M_{i_1}\cap \ldots \cap M_{i_n}$, then $I=\{i_1, \ldots,i_n\}$ is 
a positively oriented simplex in $\Gamma_M$. Furthermore, let us set 
$\sigma(p)=1\;\mathrm{or}\; -1$ according to whether the set 
$\{v_{i_1}, \ldots, v_{i_n}\}$ is a positively or negatively oriented basis of 
$H_2(BT)$ respectively. 

\begin{corollary}\label{C1}
Under the above assumptions, there exists a $T^n$-invariant almost complex structure
on the unitary torus manifold $M^{2n}$ if and only if $\sigma(p)=1$ for each $p\in M^T$.
\end{corollary}


\end{document}